\documentclass[12pt]{article}
\usepackage[english,frenchb]{babel}     
\usepackage{color}

\newtheorem{lem}{Lemme}[section]
\newtheorem{prop}{Proposition}[section]
\newtheorem{thm}{Th\'eor\`eme}[section]
\newtheorem{cor}{Corollaire}[section]{}

\newtheorem{rem}{Remarque}[section]
\newtheorem{defi}{D\'efinition}[section]

\begin{document}

\title{Borne sur le degr\'e des polyn\^omes\\ presque parfaitement non-lin\'eaires}
\author{Fran\c cois Rodier\thanks{Institut de Math\'ematiques de Luminy --
C.N.R.S. --
{e-mail:} {\tt rodier@iml.univ-mrs.fr}}}
\date{}
\maketitle


\let\fd=\rightarrow
\let\lfd=\longrightarrow
\let\lfc=\longmapsto
\let\dfd=\Rightarrow
\let\fgd=\Leftrightarrow
\let\eq=\Longleftrightarrow
\let\eps=\epsilon
\let\scr=\scriptstyle
\let\dps=\displaystyle
\let\ss=\smallskip
\let\ms=\medskip
\let\mb=\medbreak
\let\bs=\bigskip
\let\bb=\bigbreak
\let\cj=\overline

\def\Z{{\msb Z}}
\def\P{{\msb P}}
\def\f#1{{\msb F}_{#1}}
\def\mod#1{{\ ({\rm mod.}\ #1)}} 
\def\tr{\mathop{\rm Tr}\nolimits}
\def\pt#1{\left(#1\right)}
\def\ssi{\hbox{ si et seulement si }}
\def\et{\quad {\rm et\ }\quad }
\def\biq#1{\quad \hbox{ #1 }\quad }
\def\ie{c'est-\`a-dire }
\def\no{n$^\circ$}
\def\expo#1{\hbox{$^{\hbox{\tiny{#1}}}$}}

\font\tenmsb=msbm10 scaled 1200
\font\sevenmsb=msbm8
\font\fivemsb=msbm6
\newfam\msbfam
\def\msb{\fam\msbfam\tenmsb}%
\textfont\msbfam=\tenmsb \scriptfont\msbfam=\sevenmsb%
\scriptscriptfont\msbfam=\fivemsb%

\def\hrp{\ms\hrule width 2cm\mb}

\def\er#1{{\color{red}#1}}
\def\en#1{#1}

\begin{abstract}
Nyberg a d\'efini la notion de la non-lin\'earit\'e presque parfaite (APN) pour caract\'eriser les fonctions qui ont la meilleure r\'esistance aux attaques diff\'erentielles. 
Nous montrons ici que pour beaucoup de   fonctions bool\'eennes APN, le nombre de ses variables est born\'e par une expression d\'ependant du degr\'e de cette fonction. 
\end{abstract}

\selectlanguage{english}
\renewcommand{\abstractname}{English extended abstract}
\begin{abstract}
The vectorial Boolean functions are employed in cryptography to build block coding algorithms. An important criterion on these functions is their resistance to the differential cryptanalysis. Nyberg defined the notion of almost perfect non-linearity (APN) to study resistance to the differential attacks. Up to now, the study of functions APN was especially devoted to the function powers. Recently some people showed that certain quadratic polynomials were APN. 

Here, we will give a criterion so that a function is not almost perfectly non-linear. 

H. Janwa showed, by using Weil's bound, that certain cyclic codes could not correct two errors. A. Canteaut showed by using the same method that the functions powers were not APN for a too large value of the exponent. We use Lang and Weil's bound and a result of P.~Deligne (or more exactly an improvement given by Ghorpade and Lachaud) on the Weil's conjectures about surfaces on finite fields to generalize this result to many polynomials. We show therefore that these polynomials cannot be APN if their degrees are too large. We study many examples, and make some computation, showing  that one does not get any new APN function for a polynomial on a field with at most 512 elements.
\end{abstract}
\selectlanguage{francais}

\section{Introduction}

Les fonctions bool\'eennes vectorielles sont employ\'ees dans la cryptographie pour construire des algorithmes de chiffrement par bloc. Un crit\`ere important sur ces fonctions est une r\'esistance \'elev\'ee \`a la cryptanalyse diff\'erentielle. Nyberg \cite{ny} a d\'efini la notion de la non-lin\'earit\'e presque parfaite (APN) qui caract\'erise les fonctions qui ont la meilleure r\'esistance aux attaques diff\'erentielles.

 Une fonction  bool\'eenne $f$ \`a $m$ variables est d'autant plus r\'esistante aux attaques diff\'erentielles que la valeur de $\delta$ est plus petite o\`u  $$\delta=\sup_{\alpha\ne0,\beta}\#\{x\in\f2^m\mid f(x+\alpha)+f(x)=\beta \}.$$
 Les fonctions APN sont celles qui atteignent la plus petite valeur de $\delta$, \ie 2.

Jusqu'ici, l'\'etude des fonctions APN  a surtout \'et\'e consacr\'ee aux fonctions puissances (voir par exemple \cite{cpd,do1,do2,do3}).
 R\'ecemment l'\'etude s'est \'elargie \`a d'autres fonctions, en particulier aux polyn\^omes quadratiques  (Pott, Carlet et al.\cite{ekp,bcfl}) ou \`a des polyn\^omes sur des petits corps (Dillon \cite{di}).

Divers auteurs (Berger, Canteaut, Charpin,  Laigle-Chapuy \cite{bccl}, Byrne  et McGuire \cite{bg}, Jedlicka \cite{je} ou Voloch \cite{vo}) ont d\'emontr\'e l'impossibilit\'e de la propri\'et\'e pour une fonction d'\^etre APN dans certains cas.

\mb
Nous montrons ici que pour beaucoup de   fonctions polynomiales APN $f$ sur $\f {2^m}$~, le nombre $m$ est born\'e par une expression d\'ependant du degr\'e de $f$. 

Nous utilisons pour cela une m\'ethode d\'ej\`a initi\'ee par
H. Janwa qui a montr\'e, en utilisant la borne de Weil, que certains codes cycliques ne pouvaient pas corriger deux erreurs \cite{jw}. 
A. Canteaut a montr\'e en utilisant la  m\^eme m\'ethode que certaines fonctions puissances n'\'etaient pas APN pour une valeur trop grande de l'exposant \cite{ca}.
Nous avons pu g\'en\'eraliser ce r\'esultat \`a tous les polyn\^omes 
en utilisant des r\'esultats de Lang-Weil et de P.~Deligne (ou plus exactement des am\'eliorations, due \`a Ghorpade et Lachaud) sur les conjectures de Weil.

Nous terminons cette article par l'\'etude des polyn\^omes de petit degr\'e.

Une partie de ces r\'esultats a fait l'objet d'une publication partielle pour le colloque BFCA08.

\section{Pr\'eliminaires}

Pour $a\in {\f2^m}^*$, l'\'equation
$f(x+a)+f(x)=b$
a \'evidemment un nombre pair de solutions, et en a au moins une pour un $b$ dans $\f2^m$. D'o\`u la d\'efinition suivante.

\begin{defi}
Une fonction
$f:\f2^m\lfd\f2^m$
est APN 
\ssi pour tout $a\in {\f2^m}^*$ et tout $b\in \f2^m$, l'\'equation
$$f(x+a)+f(x)=b$$
n'a au plus que 2 solutions.
\end{defi}

\subsection{Polyn\^omes \'equivalents}

Posons $q=2^m$. Au lieu de fonctions dans $\f2^m$, on parlera de polyn\^omes dans $\f{q}$, \en{ce qui est \'equivalent pour les polyn\^omes de degr\'e au plus $2^m-1$.}

Un polyn\^ome $q$-affine est un polyn\^ome dont les mon\^omes sont de degr\'e 0 ou une puissance de 2.
Les propositions suivantes sont claires.

\begin{prop}
\label{p2}
La classe des fonctions APN est invariante par addition d'un polyn\^ome $q$-affine.
\end{prop}

Nous prendrons d\'esormais pour $f$ une application polynomiale  de $\f{2^m}$ dans lui m\^eme qui n'ait pas de termes de degr\'e une puissance de 2 ni de terme constant. 

\begin{prop}
\label{aff}
Pour tout $a$, $b$ et $c$ dans $\f q$ le polyn\^ome 
$cf(ax+b)$ est APN \ssi le polyn\^ome $f$ l'est.
\end{prop}

\subsection{L'\'equivalence au sens de Carlet-Charpin-Zinoviev}

Carlet, Charpin et Zinoviev ont d\'efinit une relation d'\'equivalence entre fonctions bool\'eennes  \cite{ccz}.
Pour une  fonction $f$ de $\f{2}^m$
dans lui-m\^eme on note $G_f$ le graphe de la fonction $f$:
$$G_f =\{(x, f(x)) \mid x \in\f2^m \}.$$

\begin{defi}
On dit que les fonctions $f,f' : \f2^m\lfd\f2^m$
sont \'equivalentes au sens de Carlet-
Charpin-Zinoviev (on parle aussi de
CCZ \'equivalence) s'il existe une permutation lin\'eaire $L :  \f2^{2m}\lfd\f2^{2m}$
telle que $L(G_f) = G_{f'}$ .
\end{defi}
 D'apr\`es \cite{ccz} si $f$ et $f'$ sont CCZ \'equivalentes, alors $f$ est APN si et seulement si $f'$ l'est.

\subsection{Les mon\^omes APN connus}

Il sont d\'ecrits dans la monographie de  Carlet \cite{car} o\`u on trouvera les r\'ef\'erences des d\'emonstrations.

Les fonctions $f(x)=x^d$ suivantes sont APN sur $\f{2^m}$, o\`u $d$ est donn\'e par:

\begin{itemize}
\item  $ d = 2^h + 1$ o\`u $pgcd(h, m) = 1$ (fonctions de Gold).

\item 
$d = 2^{2h}-2^h+1$ o\`u $pgcd(h, m) = 1$ (fonctions de  Kasami). 

\item 
$d = 2^{(m-1)/2} + 3$ avec $m$ impair (fonctions de Welch). 

\item  $d = 2^{(m-1)/2} + 2^{(m-1)/4} - 1$, o\`u $m\equiv 1 \mod 4$, 

  $d = 2^{(m-1)/2}+2^{(3m-1)/4}-1$, o\`u $m\equiv 3 \mod 4$ (fonctions de Niho)  

\item $d = 2^m - 2$, pour  $m$ impair;

\item 
$d = 2^{4m/5} +2^{3m/5} +2^{2m/5} +2^{m/5} - 1$, o\`u $m$ est divisible par 5 (fonctions de  Dobbertin).
\end{itemize}

Les fonctions de Gold et de Kasami sont les seules connues o\`u $d$ est ind\'ependant de $m$ et qui donnent des fonctions APN pour une infinit\'e de valeurs de $m$.

\section{Caract\'erisation des fonctions APN}

\begin{prop}
\label{apn}
La fonction
$f:\f{2^m}\lfd\f{2^m}$
est APN 
\ssi la surface affine
$$
f(x_0)+f(x_1)+f(x_2)+f(x_0+x_1+x_2)=0
$$
a tous ses points rationnels contenus dans la surface 
$(x_0+x_1)(x_2+x_1)(x_0+x_2)=0$.

\end{prop}

D\'emonstration. --

Pour que la fonction
$f:\f{2^m}\lfd\f{2^m}$
soit APN 
il faut et il suffit que
pour tout $\alpha\in\f2^m$, $\alpha\ne0$ et pour tout $\beta\in\f2^m$,
$$\#\{x_0\in\f2^m : f(x_0)+f(x_1)=\beta\ ,\ x_0+x_1=\alpha\}\le2$$
\ie que pour tout $\alpha\in\f2^m$, $\alpha\ne0$ et pour tout $\beta\in\f2^m$,
il n' y ait pas 4 \'el\'ements distincts $x_0$, $x_1$, $x_2$, $x_3$ de $\f2^m$ qui
v\'erifient
$$\left\{
\begin{array}{rcl}
x_0+x_1=\alpha, f(x_0)+f(x_1)&=&\beta\cr
x_2+x_3=\alpha, f(x_2)+f(x_3)&=&\beta
\end{array}
\right.
$$

Cela revient \`a dire qu'il n'y  a pas 3 \'el\'ements distincts $x_0$, $x_1$, $x_2$ de $\f2^m$ qui
v\'erifient
$$
f(x_0)+f(x_1)+f(x_2)+f(x_0+x_1+x_2)=0
$$
autrement dit que la surface affine
$$
f(x_0)+f(x_1)+f(x_2)+f(x_0+x_1+x_2)=0
$$
a tous ses points rationnels contenus dans la surface 
$(x_0+x_1)(x_2+x_1)(x_0+x_2)=0$.
\bs

Avant d'\'enoncer un corollaire, remarquons que
le polyn\^ome $f(x_0)+f(x_1)+f(x_2)+f(x_0+x_1+x_2)$ est divisible par 
$(x_0+x_1)(x_2+x_1)(x_0+x_2)$, par cons\'equent le quotient 
$${f(x_0)+f(x_1)+f(x_2)+f(x_0+x_1+x_2)\over (x_0+x_1)(x_2+x_1)(x_0+x_2)}$$
d\'efinit bien un polyn\^ome de degr\'e $d-3$ o\`u $d$ est le degr\'e de $f$. On v\'erifie facilement que ce polyn\^ome est nul \ssi $f$ est un polyn\^ome $q$-affine.

\begin{cor}
 Si l'application polynomiale  $f$ (de degr\'e $d\ge5$) est APN et
si  la surface affine $X$
$${f(x_0)+f(x_1)+f(x_2)+f(x_0+x_1+x_2)\over (x_0+x_1)(x_2+x_1)(x_0+x_2)}=0$$
  est absolument irr\'eductible, alors la surface projective correspondante $\cj X$ admet au plus $4((d-3)q+1)$ point rationnels, o\`u $d$ est le degr\'e de $f$ et $q=2^m$.
\end{cor}

D\'emonstration. --

Si $f$ est APN, alors $f$ n'est pas un polyn\^ome $q$-affine et l'\'equation ci-dessus d\'efinit bien une surface de degr\'e plus grand que 2.
Si la surface $\cj X$ contenait le plan $x_0+x_1=0$, elle serait \'egale \`a ce plan puisque la surface $\cj X$ est irr\'eductible et serait de degr\'e 1, ce qui est contraire \`a l'hypoth\`ese.
Par cons\'equent, elle coupe le plan $x_0+x_1=0$ suivant une courbe de degr\'e $d-3$. Cette courbe admet au plus $(d-3)q+1$ points rationnels d'apr\`es Serre \cite{se}. De m\^eme pour le plan \`a l'infini.

Si $f$ est APN, la surface $\cj X$ n'a pas d'autre points rationnels que ceux de la surface $(x_0+x_1)(x_2+x_1)(x_0+x_2)=0$, qui est r\'eunion du plan $x_0+x_1=0$  et de ses sym\'etriques, ou du plan \`a l'infini. 
Donc elle admet au plus $4((d-3)q+1)$ point rationnels.

\section{Borne inf\'erieure pour le degr\'e d'un poly\-n\^ome APN}

\subsection{Une premi\`ere borne}

\begin{thm}
\label{lawe}
Soit $f$ une application polynomiale  de $\f{2^m}$ dans lui m\^eme, $d$ son degr\'e.  Supposons que la  surface $X$ d'\'equation affine
\begin{equation}
\label{equa}
{f(x_0)+f(x_1)+f(x_2)+f(x_0+x_1+x_2)\over (x_0+x_1)(x_2+x_1)(x_0+x_2)}=0
\end{equation}
soit absolument irr\'eductible.
Alors, si $d<0,45q^{1/4}+0,5 $ et $d\ge9$\ , $f$ n'est pas APN.
\end{thm}

D\'emonstration. --

D'une am\'elioration d'un r\'esultat de Lang-Weil \cite{lw} par Ghorpade-Lachaud (\cite[section 11]{gl}), on d\'eduit
         $$|\cj X(\f{2^m})-q^2-q-1|\le (d-4)(d-5) q^{3/2}+18d^4q.$$
         D'o\`u
         $$\cj X(\f{2^m})\ge q^2+q+1-  (d-4)(d-5) q^{3/2}-18d^4q.$$
Par cons\'equent, si $ q^2+q+1-  (d-4)(d-5) q^{3/2}-18d^4q> 4((d-3)q+1)$, alors
         $\cj X(\f{2^m})>4((d-3)q+1)$, et donc $f$ n'est pas APN.
         
         Cette condition s'\'ecrit
     $$ q-(d-4)(d-5) q^{1/2}+(-18d^4-4d+13)-3/q> 0$$
         La condition est v\'erifi\'ee pour
$q^{1/2}>13,51-5 d+4,773 d^2$ si $d\ge2$.
Ou encore pour
$d<0,45q^{1/4}+0,5 $ et $d\ge9$.

\subsubsection{Irr\'eductibilit\'e de $X$}

La proposition suivante donne un crit\`ere pour que la surface $X$ soit irr\'eductible.

\begin{prop}
\label{irred}
Soit $f$ une application polynomiale  de $\f{2^m}$ dans lui m\^eme, $d$ son degr\'e.  Supposons que la courbe $X_\infty$ d'\'equation
$${x_0^d+x_1^d+x_2^d+(x_0+x_1+x_2)^d\over (x_0+x_1)(x_2+x_1)(x_0+x_2)}=0$$
soit absolument irr\'eductible.
Alors la surface $X$ d'\'equation (\ref{equa}) est absolument irr\'eductible.
\end{prop}

D\'emonstration. --

L'intersection $X_\infty$ de la surface $\cj X$ avec le plan \`a l'infini
 a comme \'equation
$${x_0^d+x_1^d+x_2^d+(x_0+x_1+x_2)^d\over (x_0+x_1)(x_2+x_1)(x_0+x_2)}=0$$
puisque $x^d$ est la composante de degr\'e $d$ du polyn\^ome $f(x)$.
Puisque la courbe $X_\infty$ est   absolument irr\'eductible il en va de m\^eme de la surface $\cj X$ donc de $X$.
\bb
Janwa, McGuire et
Wilson \cite{jgw} ont \'etudi\'e la courbe $X_\infty$ et ont d\'eduit un certain nombre de cas o\`u elle est   absolument irr\'eductible.

\begin{prop}
La courbe $X_\infty$ est absolument irr\'eductible pour les valeurs de $d\equiv 3\mod 4$ et pour les  valeurs de $d\equiv 5\mod 8$ et $d>13$.
\end{prop}

\begin{rem}
La proposition \ref{irred} donne une condition suffisante mais pas n\'ecessaire. Voir la section \ref{exe}.
\end{rem}

On peut aussi regarder l'intersection de $X$ avec le plan $x_2=a$ ou  avec le plan $x_1+x_2=a$.
\begin{prop}
Soit $f$ une application polynomiale  de $\f{2^m}$ dans lui m\^eme, $d$ son degr\'e.  Supposons que la courbe $X_a$ d'\'equation
$${f(x_0)+f(x_1)+f(a)+f(x_0+x_1+a)\over (x_0+x_1)(a+x_1)(x_0+a)}=0$$
soit absolument irr\'eductible.
Alors la surface $X$ d'\'equation (\ref{equa}) est absolument irr\'eductible.
\end{prop}

\begin{prop}
Soit $f$ une application polynomiale  de $\f{2^m}$ dans lui m\^eme, $d$ son degr\'e.  Supposons que la courbe $Y_a$ d'\'equation
$$F_a(x,y)={f(x_0)+f(x_1)+f(x_1+a)+f(x_0+a)\over (x_0+x_1)(x_0+x_1+a)}=0$$
soit absolument irr\'eductible.
Alors la surface $X$ d'\'equation (\ref{equa}) est absolument irr\'eductible.
\end{prop}
C'est la m\^eme courbe que Voloch consid\`ere dans \cite{vo}.

\subsection{Une deuxi\`eme  borne}

Sons certaines conditions, on peut obtenir une borne meilleure pour la dimension.

\begin{thm}
\label{lisse}
Soit $f$ une application polynomiale  de $\f{2^m}$ dans lui m\^eme, $d$ son degr\'e.  Supposons que la surface projective $\cj X$ associ\'ee \`a la surface $X$ d'\'equation affine
$${f(x_0)+f(x_1)+f(x_2)+f(x_0+x_1+x_2)\over (x_0+x_1)(x_2+x_1)(x_0+x_2)}=0$$
n'ait que des points singuliers isol\'es.
Alors, si $d\ge10$ et $d<q^{1/4} +4$\ , $f$ n'est pas APN.
\end{thm}

D\'emonstration. --

Notons $b_i$ le $i^{\expo{\`eme}}$ nombre de Betti de $\cj X$, \ie la dimension de l'espace de cohomologie $\ell$-adique de $\cj X$ \`a support compact, et $b'_i$ le  $i^{\expo{\`eme}}$ nombre de Betti primitif, \'egal \`a $b_i$ pour $i$ impair  et \`a $b_i-1$ pour $i$ pair. 

D'apr\`es une am\'elioration d'un r\'esultat de Deligne \cite{de} par Ghorpade-Lachaud (\cite{gl}, corollaire 7.2), on en d\'eduit
      \begin{eqnarray*}
 |X(\f{2^m})-q^2-q-1|
 &\le& b'_{1}(2, d-3) q^{3/2} + (b_2(3, d-3) + 1) q\\ 
 &\le& (d-4)(d-5) q^{3/2} + (d^3-13 d^2+57 d-82) q
\end{eqnarray*}
         D'o\`u
         $$X(\f{2^m})\ge q^2+q+1- (d-4)(d-5) q^{3/2} - (d^3-13 d^2+57 d-82) q.$$
Par cons\'equent, si 
$$q^2+q+1- (d-4)(d-5) q^{3/2} - (d^3-13 d^2+57 d-82)  q> 3((d-3)q+1),$$
 alors
         $X(\f{2^m})>4((d-3)q+1)$, et donc $f$ n'est pas APN.
        Cette condition s'\'ecrit
$$q+(-d^2+9 d-20) q^{1/2}+(-d^3+13 d^2-61 d+95) -{2\over q}> 0.$$         
       La condition est v\'erifi\'ee pour
         $ q> d^4 - 16 d^3+ 94 d^2 - 228 d +173$ d\`es que $d\ge6$.
Ou encore pour
$d<q^{1/4} +4$ d\`es que $d\ge10$.

\subsubsection{Non-singularit\'e de  $X_\infty$}

\begin{thm}
\label{surlis}
Soit $f$ une application polynomiale  de $\f{2^m}$ dans lui m\^eme, $d$ son degr\'e.  Supposons que la courbe $X_\infty$ d'\'equation
$${x_0^d+x_1^d+x_2^d+(x_0+x_1+x_2)^d\over (x_0+x_1)(x_2+x_1)(x_0+x_2)}=0$$
soit lisse.
Alors la surface $\cj X$ n'a  que des points singuliers isol\'es.
\end{thm}

{\sl Demonstration}

La d\'emonstration du th\'eor\`eme \ref{irred} montre que $X_\infty$ est l'intersection de $X$ avec l'hyperplan \`a l'infini.
Puisque la courbe $X_\infty$ est   absolument irr\'eductible il en va de m\^eme de la surface $\cj X$.

On peut d\'eduire que
 la surface $\cj X$ est r\'eguli\`ere en codimension 1 (\ie n'a  que des points singuliers isol\'es) si la courbe $X_\infty$ est non-singuli\`ere 
 (cf. Ghorpade-Lachaud \cite{gl}, 
corollaire 1.4). 

\bb
Janwa,  et
Wilson \cite{jw} ont \'etudi\'e la courbe $X_\infty$ et ont d\'eduit un certain nombre de cas o\`u elle est non singuli\`ere.

\begin{prop}
\label{lisexe}
La courbe $X_\infty$ est non singuli\`ere pour les valeurs de $d=2l+1$ o\`u 
\begin{itemize}
\item $l$ est un entier impair tel qu'il existe un entier $r$ avec $2^r\equiv-1 \mod l$.
\item $l$ est un nombre premier plus grand que 17 tel que l'ordre de 2 modulo $l$ soit $(l-1)/2$.
\end{itemize}
\end{prop}
En particulier la premi\`ere condition est satisfaite si $l$ est un nombre premier congru \`a $\pm3$ modulo 8.
Parmi les degr\'es inf\'erieurs \`a 100, on trouve 7, 11, 19, 23, 27, 35, 39, 47, 51, 55, 59, 67, 75, 83, 95.

\subsection{La surface $X$ ne peut pas \^etre lisse}

On pourrait se poser la question de savoir si la surface $\cj X$ peut \^etre lisse, ce qui am\'eliorerait encore les bornes sur le degr\'e des fonctions APN. Ce ne peut \^etre le cas.
Soit en effet
$$\phi(x_0,x_1,x_2)={f(x_0)+f(x_1)+f(x_2)+f(x_0+x_1+x_2)\over (x_0+x_1)(x_2+x_1)(x_0+x_2)}$$
l'\'equation affine de la surface $X$.
Les points singuliers de cette surface $X$ sont sur les surfaces d'\'equation
$\phi'_{x_i}(x_0,x_1,x_2)=0$ o\`u on note $\phi'_{x_i}$ la d\'eriv\'ee de $\phi$ par rapport \`a ${x_i}$.

\begin{lem}
Le polyn\^ome $x_1+x_2$ divise  $\phi'_{x_0}(x_0,x_1,x_2)=0$.
\end{lem}

{\sl D\'emonstration}

Il suffit de montrer le lemme pour chaque mon\^ome de $f$.
Si on fait le changement de variable $x_1+x_2=t$, le polyn\^ome $\phi$ devient, pour un mon\^ome de degr\'e $r$
\begin{eqnarray*}
\psi(x_0,x_1,t)&=&{x_0^r+x_1^r+(x_1+t)^r+(x_0+t)^r\over (x_0+x_1)(x_0+x_1+t)t}\\
&=&{r(x_1^{r-1}+x_0^{r-1})+tg(x_0,x_1,t)\over  (x_0+x_1)(x_0+x_1+t)}
\end{eqnarray*}
o\`u $g(x_0,x_1,t)$ est un polyn\^ome.
En d\'erivant par rapport \`a $x_0$, cela donne
\begin{eqnarray*}
&&\psi'_{x_0}(x_0,x_1,t)= {r(r-1)( x_0^{r} + x_0^{r-2} x_1^2)+tg_1(x_0,x_1,t) 
\over
(x_0 + x_1)^2 (t + x_0 + x_1)^2}
\end{eqnarray*}
o\`u $g_1(x_0,x_1,t)$ est un polyn\^ome.
Cela montre le lemme car $r(r-1)\equiv0\mod 2$.
\bb

Par cons\'equent l'intersection de  la droite $x_0=x_1=x_2$ avec la surface affine $X$ est form\'e de points singuliers de $X$.

Si la surface projective $\cj X$ ne rencontre pas la droite en des points \`a distance finie, elle ont un point commun d'ordre $d-3$ \`a l'infini. 
L'\'equation
$\phi(u,u,u)=0$ n'a pas de solution, donc $\phi(u,u,u)$ est une constante non-nulle $a_0$.
Soit $\Phi(x_0:x_1:x_2:z)$ l'\'equation projective de la surface $\cj X$. Elle est \'egale \`a $z^{d-3} \phi({x_0\over z},{x_1\over z},{x_2\over z})$, donc on a
 $\Phi(1,1,1,z)=z^{d-3} \phi({1\over z},{1\over z},{1\over z})=z^{d-3} a_0$.
La formule
$$x_0\Phi'_{x_0}+x_1\Phi'_{x_1}+x_2\Phi'_{x_2}+z\Phi'_{z}=(d-3) \Phi(x_0,x_1,x_2,z)$$
restreinte \`a la droite $x_0=x_1=x_2$
devient donc
$z\Phi'_{z}=(d-3) z^{d-3} a_0$, d'o\`u
$\Phi'_{z}=(d-3) z^{d-4} a_0$ sur cette droite.
Par cons\'equent $\Phi'_{z}(1,1,1,0)=0$ et le point $(1,1,1,0)$ est bien singulier si $d\ge5$.

\section{Autres exemples}
\label{exe}

\subsection{Bin\^omes}

\begin{prop}
\label{binome}
Soient $d$ et $r$ deux entiers tels que $d> r\ge 3$, et soit $f(x)=x^d+ax^r$, avec $a\in\f q$. Soit $\phi_s$ le polyn\^ome 
$${x_0^s+x_1^s+x_2^s+(x_0+x_1+x_2)^s\over (x_0+x_1)(x_2+x_1)(x_0+x_2)}$$
pour $s=d$ ou $r$. Supposons que $(\phi_d,\phi_r)=1$ et que
\begin{itemize}
\item ou bien $\phi_d$ se d\'ecompose en
facteurs distincts sur $\cj\f{2^m}$ et $r\ge5$;
\item ou bien  $\phi_r$ se d\'ecompose en
facteurs distincts sur  $\cj\f{2^m}$.\end{itemize}
Alors, si $d<0,45q^{1/4}+0,5 $ et $d\ge9$\ , $f$ n'est pas APN.
\end{prop}

{\sl D\'emonstration}

L'\'equation de la surface $X$ associ\'ee \`a $f$ est $\phi_d+a\phi_r=0$. Le lemme suivant montre que cette surface est irr\'eductible sous les hypoth\`eses faites.
\bb

\begin{lem}
\label{somme}
Soit $\Phi(x, y,z) \in \f{2^m}[x, y,z]$ la somme de deux
polyn\^omes homog\`enes, \ie $\Phi = \Phi_r +\Phi_d$ 
o\`u $\Phi_i$ est homog\`ene de degr\'e
$i$, et $r < d$. Supposons que $(\Phi_r, \Phi_d) = 1$ et que
\begin{itemize}
\item ou bien $\Phi_r$ se d\'ecompose en
facteurs distincts sur $\cj\f{2^m}$ et $r\ge1$;
\item ou bien  $\Phi_d$ se d\'ecompose en
facteurs distincts sur  $\cj\f{2^m}$ et $r\ge0$.\end{itemize}
Alors $\Phi$ est absolument irr\'{e}ductible sur $\f{2^m}$.
\end{lem} 

{\sl D\'emonstration} --

On fait le m\^eme raisonnement que Byrne  et McGuire dans l'article \cite[lemme 2]{bg}. La preuve fonctionne aussi dans le cas des polyn\^omes \`a trois variables. 
Elle fonctionne aussi dans le deuxi\`eme cas du lemme ($\Phi_d$ se d\'ecompose en
facteurs distincts sur  $\cj\f{2^m}$ et $r\ge0$).
\bb

Par exemple, la proposition \ref{binome} montre que le polyn\^ome
$x^{13}+a x^7$ avec $a\ne0$ ne peut \^etre APN que si  $m\le19$, car le polyn\^ome $\phi_7$ est irr\'eductible et ne divise pas $\phi_{13}$ d'apr\`es \cite{jw}.
\bb

On a aussi  la proposition suivante due pour l'essentiel \`a F. Voloch.
\begin{prop}
Soit $f(x) = x^d + cx^r$, o\`u $c \in\cj\f2^*$, $r < d$ sont des entiers, non tous les deux pairs, et non plus une puissance de 2 et tels que $(d-1, r-1)$ soit une puissance de 2. Alors, si $d<0,45q^{1/4}+0,5 $ et $d\ge9$\ , $f$ n'est pas APN.
\end{prop}

{\sl D\'emonstration} --

F. Voloch a montr\'e dans \cite{vo} que ces hypoth\`eses impliquaient que la courbe
$$F_a(x,y)={f(x_0)+f(x_1)+f(x_1+a)+f(x_0+a)\over (x_0+x_1)(x_0+x_1+a)}$$ est
irr\'{e}ductible dans $\f2[x, y,a]$.

\begin{rem}
L'\'enonc\'e de cette proposition \cite[Theorem 3]{vo} supposait que $r$ et $d$ \'etaient premiers entre eux.
F. Voloch m'a communiqu\'e que ce n'\'etait pas n\'ecessaire.
\end{rem}
\subsection{Polyn\^omes de degr\'e 3 ou 5}

D'apr\`es la proposition \ref{p2}, il suffit de regarder les polyn\^omes de la forme
$a_5 x^5 +a_3 x^3$. Ces polyn\^omes sont combinaison lin\'eaires de plusieurs mon\^omes de la forme $x^{2^i+1}$ et ils ne peuvent pas \^etre APN d'apr\`es \cite{bccl} sauf si $a_3$ ou $a_5$ est nul, auquel cas ce sont des fonctions de Gold.

\subsection{Polyn\^omes de degr\'e 6}
\label{deg6}

\begin{prop}
Soit
{$f(x)=x^6 + a_5 x^5 + a_3 x^3 $ un polyn\^ome de degr\'e 6.}
Le polyn\^ome $f$ est APN \ssi $a_3=a_5=0$.
Il est alors \'equivalent \`a un fonction de Gold.
\end{prop}

{\sl D\'emonstration} --

Il est facile de v\'erifier que la surface $X$ associ\'ee ne contient un hyperplan que si $a_3=a_5^3$.
Donc elle est absolument irr\'eductible, sauf si $a_3=a_5^3$. Dans ce cas-ci, on v\'erifie que la surface $X$ se d\'ecompose en 3 hyperplans d'\'equation
$x_0 + x_2 + a_5=0$, $ x_0 + x_1 + a_5=0$, $ a_5 + x_1 + x_2=0$
et la fonction $f$ n'est pas APN si $a_5\ne0$, d'apr\`es la proposition \ref{apn}.

Si  $a_3\ne a_5^3$ et $a_3=0$, 
on a donc
$f(x)=x^6 + a_5 x^5 $.
Comme $a_5\ne0$, on se ram\`ene par une transformation affine (cf. proposition \ref{aff}) \`a $a_5^{-6} f(a_5 x)=x^6 +  x^5  $.
Des points de la surface affine $X$ sont alors
$$\pt{{1\over\lambda (1+\lambda)},{\lambda^3 \over\lambda (1+\lambda)},1}$$
avec $\lambda\in\f q-\f2$. Si $m\ge 3$, ils ne sont pas sur la surface d'\'equation
$(x_0+x_1)(x_2+x_1)(x_0+x_2)=0,$
donc la fonction $f(x)=x^6 + a_5 x^5 $ avec $a_5\ne0$ n'est jamais APN, d'apr\`es la proposition \ref{apn}.

Si $a_3\not=a_5^3$ et si, de plus, $a_3\not=0$, alors la surface $\cj X$ n'a que des singularit\'es isol\'ees et la fonction $f$ ne peut \^etre APN que si $m\le4$ d'apr\`es la d\'emonstration du th\'eor\`eme \ref{lisse}.
D'apr\`es \cite{bl}, ces fonctions ne peuvent pas \^etre APN.

\subsection{Polyn\^omes de degr\'e 7}

\begin{prop}
Soit
$f$ un polyn\^ome de degr\'e 7.
Pour $m\ge3$, le polyn\^ome $f$ ne peut \^etre APN  que s'il est  CCZ-\'equivalent au polyn\^ome $x^7$ sur $\f{32}$.
Il est alors \'equivalent \`a un fonction de Welsh.
\end{prop}

{\sl D\'emonstration} --

Le th\'eor\`eme \ref{surlis} et la proposition \ref{lisexe} montrent que la surface $\cj X$ n'a que des points singuliers isol\'es. D'apr\`es la d\'emonstration du th\'eor\`eme \ref{lisse} on en d\'eduit que $f$ ne peut \^etre APN que si $m\le6$. La proposition s'ensuit, d'apr\`es \cite{bl} et avec une recherche exhaustive pour $m=6$.

\subsection{Polyn\^omes de degr\'e 9}

\begin{prop}
Soit
$f$ un polyn\^ome de degr\'e 9.
Le polyn\^ome $f$ ne peut \^etre APN pour une infinit\'e de $m$ que s'il est \'equivalent au polyn\^ome $x^9$.
Il est alors CCZ-\'equivalent \`a un fonction de Gold. Autrement le polyn\^ome $f$ ne peut \^etre APN que pour $m=6$ et il est \'egal \`a une fonction $f=x^9 +a_6x^6 +a_3x^3$ ou \`a une fonction CCZ-\'equivalente.
\end{prop}

{\sl D\'emonstration} --

On peut se limiter gr‰ce aux proposition \ref{p2} et \ref{aff} aux polyn\^omes 
$$f(x)=x^9+a_7x^7+a_6x^6+a_5x^5+a_3x^3.$$
Si $a_7=0$, on obtient
$$f(x)=x^9+a_6x^6+a_5x^5+a_3x^3.$$
Si, de plus, $a_6=0$, on obtient un polyn\^ome $f(x)=x^9+a_5x^5+a_3x^3$ dont tous les degr\'es des mon\^omes sont de la forme $2^r+1$. D'apr\`es \cite{bccl} il ne peut par \^etre APN, sauf si c'est un mon\^ome.
Si $a_6\ne0$,   on peut r\'eduire l'\'etude, gr‰ce \`a la proposition \ref{aff}, aux polyn\^omes
$f(x)=x^9 +a_6x^6+ x^5+a_3x^3$ et $f=x^9 +a_6x^6 +a_3x^3$.

Si $a_7\ne0$  on obtient, d'apr\`es les m\^emes propositions  le polyn\^ome
$$f(x)= x^9+ x^7+a_5x^5+a_3x^3.$$
Pour $f(x)= x^9+ x^7+a_5x^5+a_3x^3$,
on v\'erifie ais\'ement que la surface $\cj X$ n'a que des singularit\'es isol\'ees (voir ci-dessous).
La fonction $f$ ne peut \^etre APN que si $m\le13$, d'apr\`es la d\'emonstration du th\'eor\`eme \ref{lisse}.
La fonction $x^9$ est un fonction de Gold.
Pour 
 $f(x)=x^9 +a_6x^6 +a_3x^3$  avec  $a_3=a_6^2\ne0$ la fonction $f$ ne peut \^etre APN;
 pour 
 $f(x)=x^9 +a_6x^6+a_5 x^5+a_3x^3$  avec  $a_3$, $a_5$ ou $a_6$ non nuls, la fonction $f$ ne peut \^etre APN que si $m\le8$ (voir \ref{x9}). 
 Un examen exhaustif des cas restant montre la proposition.
 La seule fonction APN suppl\'ementaire que l'on trouve est une fonction $f=x^9 +a_6x^6 +a_3x^3$ pour $m=6$ d\'ej\`a obtenue par Dillon \cite{di}.
 
\subsubsection{Points singuliers de $X$ associ\'e \`a $ x^9+ x^7+a_5x^5+a_3x^3.$}

La surface $\cj X$ a pour \'equation projective
$ B_6+ B_4z^2+a_5B_2z^4+a_3z^6$
o\`u les $B_i$ sont des polyn\^omes homog\`enes de degr\'e $i$ en $x_0,x_1,x_2$ et ce sont les \'equations des courbes $X_\infty$ correspondantes aux degr\'es $i$.
Par \cite{jw}, les polyn\^omes $B_6$ et $B_2$ sont des \'equations de r\'eunions de droites passant par le point $(1:1:1)$, d\'efinies sur $\f8$ et $\f4$ et le polyn\^ome $B_4$ est l'\'equation d'une courbe irr\'eductible.
D'apr\`es le calcul des d\'eriv\'ees de $X$, on \'etudie plusieurs cas.

Cas \no\ 1:
$x_0=x_1=x_2$.
Les points de la surface $X$ qui sont sur cette droite sont ceux qui v\'erifient
$ x_0^4+a_3z^4=0$ ou $z=0$.

Cas \no\ 2:
$x_1=x_2$.
Les points singuliers de la surface $X$ qui sont sur ce plan sont ceux qui v\'erifient
$x_0^4+x_2^4+x_0 x_2 z^2+a_5 z^4=0$
et
$x_0^6 + x_0^4 x_2^2 + x_0^2 x_2^4 + x_2^6 + x_0^4 z^2 + x_0^2 x_2^2 z^2 + 
 x_2^4 z^2 + a_5 x_0^2 z^4 + a_5 x_2^2 z^4 + a_3 z^6=0$
 qui est 2 fois la courbe
$(x_0+ x_2)^3 + (x_0^2   + x_0  x_2    + 
 x_2^2 )z  + a_5( x_0  + x_2 ) z^2+ a_3 z^3=0$.
 Un changement de variable $s=x_0+ x_2$ donne, pour $z=1$:
 $s^4 +x_0( x_0+s)  +a_5  =0$
 et
 $s^3 +  s^2   + x_0  (s+x_0) + a_5 s+ a_3 =0$  qui sont deux courbes elliptiques distinctes.
 Les deux courbes ne se coupent qu'en un nombre fini de points.

Cas \no\ 3:
$x_0\ne x_1\ne x_2\ne x_0$.

Les points singuliers de la surface $X$ sont ceux qui v\'erifient
$$\left\{
\begin{array}{ccc}
a_5+x_1 x_2+B_2^2(x_0,x_1,x_2) &=& 0\\
a_5+x_0 x_2+B_2^2(x_0,x_1,x_2) &=& 0\\
a_5+x_1 x_0+B_2^2(x_0,x_1,x_2) &=& 0
\end{array}
\right.$$
donc, par soustraction, qui v\'erifient
$$\left\{
\begin{array}{ccc}
(x_1 +x_0) x_2&=&0\\
(x_2+x_1) x_0&=&0\\
(x_0+ x_2)x_1&=&0\\
\end{array}
\right.$$
D'apr\`es l'hypoth\`ese, on a donc
$x_0= x_1= x_2$, donc on n'a rien de mieux.
\subsubsection{Surface associ\'ee \`a $ x^9+a_6x^6+a_5x_5+a_3x^3$.}
\label{x9}

Si $a_3=a_5=a_6=0$, la fonction $x^9$ est une fonction de Gold.
Supposons que l'un de ces coefficient soit non nul.

La surface $X$ a pour \'equation projective
$ B_6+ a_6B_3z^3+ a_5B_2z^4+a_3z^6$
o\`u les $B_i$ sont comme plus haut. 
L'intersection de la surface $X$ avec le plan $x_2=0$ est la courbe $Y$ d'\'equation affine
$$\phi(x_0,x_1,0)= B_6(x_0,x_1,0)+ a_6B_3(x_0,x_1,0)+ a_5B_2(x_0,x_1,0)+a_3$$
avec
$B_6(x_0,x_1,0)=\prod_{\beta\in\f8-\f4}(x_0+\beta x_1)$,
$B_3(x_0,x_1,0)=x_0x_1(x_0+x_1)$ et
$B_2(x_0,x_1,0)=(x_0+\alpha x_1)(x_0+\alpha^2x_1)$ o\`u $\alpha\in\f4-\f2$.

Les raisonnements de Byrne et McGuire pour prouver le  lemme 2 dans \cite{bg} peuvent s'appliquer.
On trouve ainsi  que le seul cas o\`u $\phi(x_0,x_1,0)$ pourrait \^etre r\'eductible est:
$$\phi(x_0,x_1,0)= \pt{ P_3+P_0} \pt{Q_3+Q_0} =P_3 Q_3+P_3Q_0 +P_0Q_3+P_0Q_0$$
o\`u les $P_i$ et les $Q_i$ soient des polyn\^omes homog\`enes de degr\'e $i$. De plus le fait que $P_3Q_0 +P_0Q_3= a_6B_3(x_0,x_1,0)$ implique que 
$P_3=(x_0+\beta x_1)(x_0+\beta^2 x_1)(x_0+\beta^4 x_1)$,
$Q_3=(x_0+\beta^3 x_1)(x_0+\beta^5 x_1)(x_0+\beta^6 x_1)$,
et que $P_0 =Q_0$.
Cela implique d\'ej\`a que si $a_5\ne0$, la courbe $Y$ est ind\'ecomposable.

On trouve donc que le seul cas o\`u la courbe $Y$  est r\'eductible est le cas o\`u $a_3=a_6^2$ et $a_5=0$.
Dans ce cas, la surface $X$ se d\'ecompose en $X=X_1\cup X_2$, o\`u $X_1$ a pour \'equation
$x_0^3 + x_0^2 x_1 + x_1^3 + x_1^2 x_2 + x_0 x_2^2 + x_2^3+a_6 $
et
$X_2$ a pour \'equation
$x_0^3 + x_0 x_1^2 + x_1^3 + x_0^2 x_2 + x_1 x_2^2 + x_2^3+a_6=0$.
L'intersection de la surface $X_1$ avec le plan $x_2=0$ est une courbe elliptique d'\'equation affine
$1 + x_0^3 + x_0 x_1^2 + x_1^3=0$.
Elle ne peut rencontrer les droites $x_1=0$, $x_2=0$, $x_1+x_3=0$ ou la droite \`a l'infini qu'en 3 points chacune. Or elle a au moins
$1+q-2\sqrt q$ points rationnels d'apr\`es la borne de Hasse-Weil, ce qui prouve que la surface affine $X_1$ a des points rationnels en dehors de la surface 
$(x_0+x_1)(x_2+x_1)(x_0+x_2)=0$ pour $1+q-2\sqrt q>12$.
Par cons\'equent, $f(x)= x^9+a_6x^6+a_6^2x^3$ n'est pas APN pour $m\ge5$. Pour $m\le4$, le calcul exhaustif montre que l'on obtient pas de nouvelles fonctions APN (cf. \cite{bl}).

Dans le cas o\`u la courbe $Y$ est irr\'eductible, elle a au plus
$1+q-20\sqrt q$ points rationnels d'apr\`es le corollaire 7.4 de \cite{gl}.
Si la fonction $f$ est APN, la courbe projective $\cj Y$ doit avoir tous ses points rationnels contenus dans la courbe r\'eunion de 4 droites $x_0x_1(x_0+x_1)z=0$.
Elle recoupe chaque droite en 6 points au plus, donc elle doit avoir au plus 24 points. Cela n'est possible que si $m\le 8$.

\section{Applications num\'eriques}

Quand la surface $X$ est irr\'eductible, on obtient que la fonction $f$ ne peut \^etre APN que si $m\le m_{max}$ o\`u $m_{max}$ est donn\'e par le tableau suivant.

\begin{center}
\begin{tabular}{|c|c|c|c|c|c|c|c|c|c|c|c|c|c|c|c|c|c|c|c|}
\hline
$d\le$&7&9&10&12&15&17&21&23&29& 36& 41&  49&50&70&83\\
\hline
$m_{max}$&15&16&17&18&19&20&21&22&23&24& 25& 26& 27&28&29\\
\hline
\end{tabular}
\end{center}

Quand la surface $\cj X$ est de plus \`a singularit\'es isol\'ees, on obtient que la fonction $f$ ne peut \^etre APN que si $m\le m_{max}$ o\`u $m_{max}$ est donn\'e par le tableau suivant.

\begin{center}
\begin{tabular}{|c|c|c|c|c|c|c|c|c|c|c|c|c|c|c|c|}
\hline
$d\le$&7&9&10& 12&13&15&17& 20& 23&26& 30& 36& 42&49&57\\
\hline
$m_{max}$&6&9&10& 11&12&13&14& 15& 16&17& 18& 19& 20&21&22\\
\hline
\end{tabular}
\end{center}

\section{Remerciements}

Je remercie Felipe Voloch qui m'a aid\'e par ses commentaires sur une ancienne version de l'article et Gregor Leander qui m'a aid\'e pour les recherches exhaustives.



\begin{thebibliography}{XXXX}

\bibitem{bccl}{T. Berger, A. Canteaut, P. Charpin,  Y. Laigle-Chapuy}
{\it On almost perfect nonlinear functions over $F\sp n\sb 2$.}
IEEE Trans. Inform. Theory 52 (2006), no. 9, 4160--4170. 


\bibitem{bcfl}
{L. Budaghyan and C. Carlet and P. Felke and G. Leander}
{\sl An infinite class of quadratic APN functions which are not equivalent to power mappings,}
Cryptology ePrint Archive, \no\  2005/359

\bibitem{bg}
{Byrne E. and McGuire G.},
{\sl On the Non-Existence of Quadratic APN and
Crooked Functions on Finite Fields}, prepublication.
http://www.maths.may.ie/sta/gmg/APNniceWeilEBGMG.pdf.


\bibitem{bl}
{Marcus Brinkman, G. Leander}: {\sl On the classification of APN functions up to
dimension five}, International Workshop on Coding and Cryptography
(WCC), Versailles, France, 2007.

\bibitem{ca}
{Canteaut, A.}, {\sl Differential cryptanalysis of Feistel ciphers
and differentially $\delta$-uniform mappings}.
 In Selected Areas on Cryptography,
     SAC'97, pages 172-184, Ottawa, Canada, 1997.
 
\bibitem{cpd}
 {A. Canteaut, P. Charpin, and H. Dobbertin.} {\sl Weight divisibility of
     cyclic codes, highly nonlinear functions on $GF(2^m)$ and crosscorrelation
     of maximum-length sequences}. SIAM Journal on Discrete Mathematics,
     13(1), 2000.
     
\bibitem{car}
     {C. Carlet.} {\sl Vectorial Boolean Functions for Cryptography.} Chapter of the monography
Boolean Methods and Models, Y. Crama and P. Hammer eds, Cambridge University
Press, to appear.

\bibitem{ccz}
{C. Carlet, P. Charpin and V. Zinoviev}. {\sl Codes, bent functions and permutations
suitable for DES-like cryptosystems.} Designs, Codes and Cryptography, 15(2), pp.
125-156, 1998.

\bibitem{de}
{Deligne, Pierre} 
{\sl La conjecture de Weil : I.} Publications Math\'ema\-tiques de l'IHES, 43 (1974), p. 273-307 

\bibitem{di}
{J. F. Dillon}
{\sl APN Polynomials and Related Codes}
Conference on Polynomials over Finite Fields
and Applications,
Banff International Research Station
November 2006.

\bibitem{do1}
{H. Dobbertin}. {\sl Almost perfect nonlinear power functions over $GF(2^n)$: the Niho case.}
Inform. and Comput., 151, pp. 57-72, 1999.

\bibitem{do2}
{H. Dobbertin}. {\sl Almost perfect nonlinear power functions over $GF(2^n)$: the Welch case.}
IEEE Trans. Inform. Theory, 45, pp. 1271-1275, 1999.

\bibitem{do3}
{H. Dobbertin}. {\sl Almost perfect nonlinear power functions over $GF(2^n)$: a new case for
$n$ divisible by 5.} D. Jungnickel and H. Niederreiter eds. Proceedings of Finite Fields
and Applications FQ5, Augsburg, Germany, Springer, pp. 113-121, 2000.

\bibitem{ekp}
{Y. Edel, G. Kyureghyan and A. Pott.} 
{\sl A new APN function which is not equivalent
to a power mapping.} Preprint, 2005, http://arxiv.org/abs/math.CO/0506420

\bibitem{gl}
{Ghorpade Sudhir R., Lachaud Gilles.}
         {\sl \'Etale cohomology, Lefschetz theorems and number of points of singular          varieties over finite fields.}, 
           Mosc. Math. J.  2  (2002),  \no\  3, 589--631.
           
\bibitem{je}
{Jedlicka, D.},
{\sl APN monomials over ${\rm GF}(2\sp n)$ for infinitely many $n$},
Finite Fields Appl. 13 (2007), no. 4, 1006--1028. 
         
\bibitem{jgw}
           {Janwa, H., McGuire, G., Wilson, R.},
            {\sl Double-error-correcting cyclic codes and absolutely irreducible polynomials over ${\rm GF}(2)$}. J. Algebra 178 (1995), no. 2, 665--676.
           
           
\bibitem{jw}
         {Janwa, H.,
Wilson, R. M}. 
{\sl Hyperplane sections of Fermat varieties in $P^3$ in char. 2
and some applications to cyclic codes.} 
         Applied algebra, algebraic algorithms and error-correcting codes
(San Juan, PR, 1993), 180--194, Lecture Notes in Comput. Sci., 673,
Springer, Berlin, 1993.


         
\bibitem{lw}
         {Lang, Serge, Weil, Andr\'e}
{\sl Number of points of varieties in finite fields}, 
Amer. J. Math. 76, (1954). 819--827.

         \bibitem{ny}
{Nyberg, Kaisa}, {\sl Differentially uniform mappings for cryptography.} Advances in
 cryptology---EUROCRYPT '93 (Lofthus, 1993), 55--64, Lecture Notes in Comput. Sci., 765, Springer, Berlin,
 1994.

         \bibitem{se}
{J. -P. Serre}, {\sl Lettre  \`a M. Tsfasman}, Ast\'erisque 198-199-200 (1991), 351-353.

         \bibitem{vo}
{Voloch, F.},
{\sl Symmetric Cryptography and Algebraic Curves,} preprint,
http://www.ma.utexas.edu/users/voloch/preprint.html

\end{thebibliography}
         \end{document}